\magnification 1200
\def\O{\Omega}
\def\e{\epsilon}
\def\w{\wedge}

\input amstex

\centerline{\bf Some properties of the complex Monge-Amp\`{e}re operator in}
\centerline{\bf Cegrell's classes and applications}
\vskip0,5cm
\noindent
\centerline{NGUYEN VAN KHUE and PHAM HOANG HIEP}
\vskip0,5cm
\noindent
{\bf Abstract.} In this article we will first prove a result about convergence in capacity. Using the achieved result we will obtain a general decompositon theorem for complex Monge-Amp\`{e}re measues which will be used to prove a comparison principle for the complex Monge-Amp\`{e}re operator.
\medskip
\noindent
2000 Mathematics Subject Classification: Primary 32W20, Secondary 32U15.
\medskip
\noindent
Key words and phrases: complex Monge-Amp\`{e}re operator, plurisubharmonic function.
\medskip
\noindent
This work was supported by the National Research Program for Natural Sciences, Vietnam.
\medskip
\noindent
{\bf 1. Introduction}
\medskip
\noindent
Let $\O$ be a bounded hyperconvex domain in $\bold C^n$. By PSH$(\O)$ we denote the set of plurisubharmonic (psh) functions on $\O$. In [BT 1,2] the authors established and used the comparison principle to study the Dirichlet problem in $\text{PSH}\cap L_{loc}^{\infty}(\O)$. Recently, Cegrell introduced a general class $\Cal E$ of psh functions on which the complex Monge-Amp\`{e}re operator $(dd^c.)^n$ can be defined. He obtained many important results of pluripotential theory in the class $\Cal E$. For example, the ones on the comparison principle and solvability of the Dirichlet problem (see [Ce 1-3]).
\medskip
\noindent
The main result of our paper are Theorem 4.1 and some Xing type comparision principles. Theorem 4.1 is generalize Lemma 5.4 in [Ce1], Lemma 7.2 in [\AA h] and Lemma 3.4 in [Ce3]. For definitions of Cegrell's classes see Section 2. After giving some preliminaries, we start in Proposition 3.1 with a comparison principle, which is analogous to a comparison principle due to Xing (Lemma 1 in [Xi1]). It should be observed that our proof is quite different from Xing's proof, and the inequality we obtain is slightly stronger than Xing's inequality, even in the case of bounded psh functions. Using Proposition 3.1, we give in Theorem 3.5 a sufficient condition for $C_n$-capacity convergence of a sequence of psh functions in the class $\Cal F$. This result should be compared to Theorem 3 of [Xi1] where the situation of bounded psh functions was studied. Applying Theorem 3.5 we give generalizations of recent results in [Cz] and [CLP] about convergences of multipole Green functions and a criterion for pluripolarity, respectively. Section 4 focuses on Theorem 4.1 and Theorem 4.9. By applying Theorem 4.1 we give some results on class Cegrell's classes. We prove in Proposition 4.4 a local estimate for the Monge-Amp\`{e}re measure in terms of the Beford-Taylor relative capacity. As an application, we give in Theorem 4.5 a decomposition result for Monge-Amp\`{e}re measure, which is similar in spirit to Theorem 6.3 in [Ce1]. From  Proposition 3.1 and Theorem 4.1 we obtain easily a Xing type comparison principle for functions in classes $\Cal F$ and $\Cal E$.
\medskip
\noindent
{\bf Acknowledgment.} We are grateful to Professor Urban Cegrell for useful  discussions that helped to improve the paper. We are grateful to Per \AA hag for fruitful comments. This work is supported by the National Research Program for Natural Sciences, Vietnam.
\vskip1cm
\noindent
{\bf 2. Preliminaries}
\medskip
\noindent
First we recall some elements of pluripotential theory that will be used throughout the paper. All this can be found in [BT2], [Ce1], [Ce2], [Le].
\medskip
\noindent
{\bf 2.1.} We will always denote by $\O$ a bounded hyperconvex domain in $\bold C^n$ unless other wise stated. The $C_n$-capacity in the sense of Bedford and Taylor on $\O$ is the set function given by
$$C_n(E)=C_n(E,\O)=\sup\{\int\limits_{E}(dd^cu)^n\ :\ u\in\text{PSH}(\O),\ -1\leq u\leq 0\}$$
for every Borel set $E$ in $\O$. It is proved in [BT2] that
$$C_n(E)=\int\limits_E (dd^ch_{E,\O}^*)^n$$
where $h_{E,\O}^*$ is the upper regularization of the relative extremal function $h_{E,\O}$ for $E$ (relative to $\O$) i.e.,
$$h_{E,\O}(z)=\sup\{u(z)\ :\ u\in\text{PSH}^-(\O),\ u\leq -1\ \text{on}\ E\}.$$
The following concepts are taken from [Xi1] and [Xi2]
\medskip
\noindent
$*$A sequence of functions $u_j$ on $\Omega$ is said to converge to a function $u$ in $C_n$-capacity on a set $E\subset\Omega$ if for every $\delta > 0$  we have $C_n(\{z\in E\ :\ |u_j(z)-u(z)|>\delta\})\to 0$ as $j\to\infty.$
\medskip
\noindent
$*$A family of positive measures \{$\mu_\alpha$\} on $\Omega$ is called uniformly absolutely continuous with respect to $C_n$-capacity in a set $E\subset\O$ if for every $\epsilon >0$ there exists $\delta >0$ such that for each Borel subset $F\subset E$ with $C_n$(F)$<\delta$ the inequality $\mu_\alpha$(F)$<\epsilon$ holds for all $\alpha$. We write $\mu_\alpha\ll C_n$ in $E$ uniformly for $\alpha$.
\medskip
\noindent
{\bf 2.2.} The following classes of psh functions were introduced by Cegrell in [Ce1] and [Ce2]
$$\Cal E_0=\Cal E_0(\Omega)=\{\varphi\in \text{PSH}^-(\Omega)\cap L^\infty(\Omega):\ \lim\limits_{z\to\partial\Omega}\varphi (z)=0,\ \int\limits_\Omega (dd^c\varphi)^n<+\infty\},$$
$$\Cal F=\Cal F(\Omega)=\{\varphi\in \text{PSH}^-(\Omega):\ \exists\ \Cal E_0(\O)\ni\varphi_j\searrow\varphi,\ \sup\limits_{j\geq 1}\int\limits_\Omega(dd^c\varphi_j)^n<+\infty\},$$
$$\Cal E=\Cal E(\Omega)=\{\varphi\in \text{PSH}^-(\Omega):\ \exists\ \varphi_K\in\Cal F(\O)\ \text{such that}\ \varphi_K=\varphi\ \text{on}\ K,\ \forall K\subset\subset\Omega\},$$
$$\Cal E^a=\Cal E^a(\O)=\{u\in\Cal E(\O)\ :\ (dd^cu)^n(E)=0\ \forall\ E\ \text{is pluripolar in}\ \O\}.$$
For each $u\in\Cal F(\O)$, we set
$$e_0(u)=\int\limits_\O (dd^cu)^n.$$
{\bf 2.3.} Let $A=\{(w_j,\nu_j)\}_{j=1,...,p}$ be a finite subset of $\O\times\bold R^+$. According to Lelong (see [Le]), the pluricomplex Green function with poles in $A$ is defined by
$$g(A)(z)=\sup\{u(z):\ u\in\Cal L_A\}$$
where
$$\Cal L_A=\{u\in\text{PSH}^-(\O):\ u(z)-\nu_j\log|z-w_j|\leq O(1)\ \text{as}\ z\to w_j,\ j=1,...,p\}$$
Set
$$\nu (A)=\sum\limits_{j=1}^p\nu_j^n,\ \hat{A}=\{w_j\}_{j=1,...,p}.$$
{\bf 2.4.} We write $\varliminf\limits_{z\to\partial\O}[u(z)-v(z)]\geq a$ if for every $\epsilon >0$ there exists a compact set K in $\O$ such that
$$u(z)-v(z)\geq a-\epsilon\ \text{for}\ z\in(\O\backslash K)\cap\{u>-\infty\}$$
and
$$v(z)=-\infty\ \text{for}\ z\in(\O\backslash K)\cap\{u=-\infty\}.$$
{\bf 2.5.} Xing's comparison principle (see Lemma 1 in [Xi1]). {\sl Let $\O$ be a bounded open subset in $\bold C^n$ and} $u,v\in\text{PSH}\cap L^\infty (\O)$ {\sl satisfy $\varliminf\limits_{z\to\partial\O}[u(z)-v(z)]\geq 0$. Then for any constant $r\geq 1$ and all} $w_j\in\text{PSH} (\O)$ {\sl with $0\leq w_j\leq 1$, $j=1,2,...,n$ we have}
$$\frac 1 {(n!)^2} \int\limits_{\{u<v\}} (v-u)^ndd^cw_1\w ...\w dd^cw_n+\int\limits_{\{u<v\}}(r-w_1)(dd^cv)^n\leq\int\limits_{\{u<v\}}(r-w_1)(dd^cu)^n$$
\vskip1cm
\noindent
{\bf 3. Some convergence theorems}
\medskip
\noindent
In order to study the convergence of a sequence of psh functions in $C_n$-capacity, we start with the following.
\medskip
\noindent
{\bf 3.1. Proposition.} a) {\sl Let $u,v\in\Cal F$such that $u\leq v$ on $\O$. Then for $1\leq k\leq n$}
$$\frac 1 {k!} \int\limits_\O(v-u)^kdd^cw_1\w ...\w dd^cw_n+\int\limits_\O(r-w_1)(dd^cv)^k\w dd^cw_{k+1}\w ...\w dd^cw_n$$
$$\leq \int\limits_\O(r-w_1)(dd^cu)^k\w dd^cw_{k+1}\w ...\w dd^cw_n$$
{\sl for all} $w_j\in\text{PSH}(\O),\ 0\leq w_j\leq 1,\ j=1,...,k,\ w_{k+1},...,w_n\in\Cal F$ {\sl and all $r\geq 1$.}
\medskip
\noindent
b) {\sl Let $u,v\in\Cal E$ such that $u\leq v$ on $\O$ and $u=v$ on $\O\backslash K$ for some $K\subset\subset\O$. Then for $1\leq k\leq n$}
$$\frac 1 {k!} \int\limits_\O(v-u)^kdd^cw_1\w ...\w dd^cw_n+\int\limits_\O(r-w_1)(dd^cv)^k\w dd^cw_{k+1}\w ...\w dd^cw_n$$
$$\leq \int\limits_\O(r-w_1)(dd^cu)^k\w dd^cw_{k+1}\w ...\w dd^cw_n$$
{\sl for all} $w_j\in\text{PSH}(\O),\ 0\leq w_j\leq 1,\ j=1,...,k,\ w_{k+1},...,w_n\in\Cal E$ {\sl and all $r\geq 1$.}
\medskip
\noindent
We proceed through some lemmas.
\medskip
\noindent
{\bf 3.2. Lemma.} {\sl Let} $u,v\in\text{PSH}\cap L^\infty(\O)$ {\sl such that $u\leq v$ on $\O$ and $\lim\limits_{z\to\partial\O}[u(z)-v(z)]=0$. Then}
\medskip
\noindent
$$\int\limits_\O (v-u)^kdd^cw\w T\leq k\int\limits_\O (1-w)(v-u)^{k-1}dd^cu\w T$$
{\sl for all} {\sl $w\in\text{PSH}(\O)$, $0\leq w\leq 1$ and all positive closed currents $T$.}
\medskip
\noindent
{\sl Proof.} First, assume $u,v\in\text{PSH}\cap L^\infty(\O)$, $u\leq v$ on $\O$ and $u=v$ on $\O\backslash K$, $K\subset\subset\O$. Then, using the Stokes formula we obtain
$$\aligned\int\limits_\O (v-u)^kdd^cw\w T&=\int\limits_\O (v-u)^kdd^c(w-1)\w T
\\&=\int\limits_\O (w-1)dd^c(v-u)^k\w T
\\&=-k(k-1)\int\limits_\O (1-w)d(v-u)\w d^c(v-u)\w T
\\&\ \ \ \ +k\int\limits_\O (1-w)(v-u)^{k-1}dd^c(u-v)\w T
\\&\leq k\int\limits_\O (1-w)(v-u)^{k-1}dd^c(u-v)\w T
\\&\leq k\int\limits_\O (1-w)(v-u)^{k-1}dd^cu\w T.\endaligned$$
General case, for each $\e>0$ we set $v_\e=\max (u,v-\e)$. Then $v_\e\nearrow v$ on $\O$, $v_\e\geq u$ on $\O$ and $v_\e =u$ on $\O\backslash K$ for some $K\subset\subset\O$. Hence
\medskip
\noindent
$$\int\limits_\O (v_\e -u)^kdd^cw\w T\leq k\int\limits_\O (1-w)(v_\e -u)^{k-1}dd^cu\w T.$$
Since $0\leq v_\e -u\nearrow v-u$ as $\e\searrow 0$, letting $\e\searrow 0$ we get
$$\int\limits_\O (v-u)^kdd^cw\w T\leq k\int\limits_\O (1-w)(v-u)^{k-1}dd^cu\w T.$$
{\bf 3.3. Lemma.} {\sl Let} $u,v\in\text{PSH}\cap L^\infty(\O)$ {\sl such that $u\leq v$ on $\O$ and $\lim\limits_{z\to\partial\O}[u(z)-v(z)]=0$. Then for $1\leq k\leq n$}
$$\frac 1 {k!} \int\limits_\O (v-u)^kdd^cw_1\w ...\w dd^cw_n+\int\limits_\O (r-w_1)(dd^cv)^k\w T$$
$$\leq \int\limits_\O (r-w_1)(dd^cu)^k\w T.$$
{\sl for all} $w_1,...,w_k\in\text{PSH}(\O),\ 0\leq w_j\leq 1\ \forall\ j=1,...,k,\ w_{k+1},...,w_n\in\Cal E$ {\sl and all $r\geq 1$.}
\medskip
\noindent
{\sl Proof.} To simplify the notation we set
$$T=dd^cw_{k+1}\w ...\w dd^cw_n.$$
First, assume that $u,v\in\text{PSH}\cap L^\infty(\O),\ u\leq v$ on $\O$, and $u=v$ on $\O\backslash K$, $K\subset\subset\O$. Using Lemma 3.2 we get
$$\aligned\int\limits_\O (v-u)^kdd^cw_1\w ...\w dd^cw_n&\leq k\int\limits_\O(v-u)^{k-1}dd^cw_1\w ...\w dd^cw_{k-1}\w dd^cu\w T
\\&\leq ...
\\&\leq k!\int\limits_\O(v-u)dd^cw_1\w (dd^cu)^{k-1}\w T
\\&\leq k!\int\limits_\O(v-u)dd^cw_1\w [\sum\limits_{i=0}^{k-1}(dd^cu)^i\w (dd^cv)^{k-i-1}]\w T
\\&= k!\int\limits_\O (w_1-r)dd^c(v-u)\w [\sum\limits_{i=0}^{k-1}(dd^cu)^i\w (dd^cv)^{k-i-1}]\w T
\\&= k!\int\limits_\O (r-w_1)dd^c(u-v)\w [\sum\limits_{i=0}^{k-1}(dd^cu)^i\w (dd^cv)^{k-i-1}]\w T
\\&= k!\int\limits_\O (r-w_1)[(dd^cu)^k-(dd^cv)^k]\w T.\endaligned$$
General case, for each $\e>0$ we put $v_\e=\max (u,v-\e)$. Then $v_\e\nearrow v$ on $\O$, $v_\e\geq u$ on $\O$ and $v_\e =u$ on $\O\backslash K$ for some $K\subset\subset\O$. Hence
$$\frac 1 {k!} \int\limits_\O(v_\e -u)^kdd^cw_1\w ...\w dd^cw_n+\int\limits_\O(r-w_1)(dd^cv_\e )^k\w T$$
$$\leq \int\limits_\O(r-w_1)(dd^cu)^k\w T.$$
Observe that $0\leq v_\e-u\nearrow v-u$ and $(dd^cv_\e)^k\w T\to (dd^cv)^k\w T$ weakly as $\e\searrow 0$, $r-w_1$ is lower semicontinuous, by letting $\e\searrow 0$ we have
$$\frac 1 {k!} \int\limits_\O (v-u)^kdd^cw_1\w ...\w dd^cw_n+\int\limits_\O (r-w_1)(dd^cv)^k\w T$$
$$\leq \int\limits_\O (r-w_1)(dd^cu)^k\w T.$$
The proof is finished.
\medskip
\noindent
{\sl Proof of Proposition 3.1.} a) Let $\Cal E_0\ni u_j\searrow u$ and $\Cal E_0\ni v_j\searrow v$ as in the definition of $\Cal F$. Replace $v_j$ by $\max(u_j,v_j)$ we may assume that $u_j\leq v_j$ for $j\geq 1$. By Lemma 3.3 we have
$$\frac 1 {k!} \int\limits_\O(v_j-u_t)^kdd^cw_1\w ...\w dd^cw_n+\int\limits_\O(r-w_1)(dd^cv_j)^k\w dd^cw_{k+1}\w ...\w dd^cw_n$$
$$\leq \int\limits_\O(r-w_1)(dd^cu_t)^k\w dd^cw_{k+1}\w ...\w dd^cw_n$$
for $t\geq j\geq 1$. By Proposition 5.1 in [Ce2] letting $t\to\infty$ in the above inequality we have
$$\frac 1 {k!} \int\limits_\O(v_j-u)^kdd^cw_1\w ...\w dd^cw_n+\int\limits_\O(r-w_1)(dd^cv_j)^k\w T$$
$$\leq \int\limits_\O(r-w_1)(dd^cu)^k\w T$$
for $j\geq 1$. Next letting $j\to\infty$ again by Proposition 5.1 in [Ce2] we get the desired conclusion.
\medskip
\noindent
b) Let $G,W$ be open sets such that $K\subset\subset G\subset\subset W\subset\subset\O$. According to the remark following Definition 4.6 in [Ce2] we can choose a function $\tilde v\in\Cal F$ such that $\tilde v\geq v$ and $\tilde v=v$ on $W$. Set
$$\tilde u=\cases & u\ \text{on}\ G
\\&\tilde v\ \text{on}\ \O\backslash G\endcases$$
Since $u=v=\tilde v$ on $W\backslash K$ we have $\tilde u\in\text{PSH}^-(\O)$. It is easy to see that $\tilde u\in\Cal F$, $\tilde u\leq\tilde v$ and $\tilde u=u$ on $W$. By a) we have
$$\frac 1 {k!} \int\limits_\O(\tilde v-\tilde u)^kdd^cw_1\w ...\w dd^cw_n+\int\limits_\O(r-w_1)(dd^c\tilde v)^k\w dd^cw_{k+1}\w ...\w dd^cw_n$$
$$\leq \int\limits_\O(r-w_1)(dd^c\tilde u)^k\w dd^cw_{k+1}\w ...\w dd^cw_n.$$
Since $\tilde u=\tilde v$ on $\O\backslash G$ we have
$$\frac 1 {k!} \int\limits_W(\tilde v-\tilde u)^kdd^cw_1\w ...\w dd^cw_n+\int\limits_W(r-w_1)(dd^c\tilde v)^k\w dd^cw_{k+1}\w ...\w dd^cw_n$$
$$\leq \int\limits_W(r-w_1)(dd^c\tilde u)^k\w dd^cw_{k+1}\w ...\w dd^cw_n.$$
Since $\tilde u=u$, $\tilde v=v$ on $W$ and $u=v$ on $\O\backslash K$ we obtain
$$\frac 1 {k!} \int\limits_\O(v-u)^kdd^cw_1\w ...\w dd^cw_n+\int\limits_\O(r-w_1)(dd^cv)^k\w dd^cw_{k+1}\w ...\w dd^cw_n$$
$$\leq \int\limits_\O(r-w_1)(dd^cu)^k\w dd^cw_{k+1}\w ...\w dd^cw_n.$$
{\bf 3.4. Proposition.} {\sl Let $u,v\in\Cal F$ and $u\leq v$ on $\O$. Then}
$$\frac 1 {n!} \int\limits_\O (v-u)^ndd^cw_1\w ...\w dd^cw_n\leq\int\limits_\O(-w_1)[(dd^cu)^n-(dd^cv)^n]$$
{\sl for all} $w_j\in\text{PSH}(\O)$, $-1\leq w_j\leq 0$, $j=1,...,n$.
\medskip
\noindent
{\sl Proof.} The proposition follows from Proposition 3.1 with $k=n,\ r=1$ and $w_j$ are replaced by $w_j+1$.
\medskip
\noindent
{\bf 3.5. Theorem.} {\sl Let $u,u_j\in\Cal F$ and $u_j\leq u$ for $j\geq 1$. Assume that $\sup\limits_{j\geq 1}\int\limits_\O(dd^cu_j)^n<+\infty$ and $||(dd^cu_j)^n-(dd^cu)^n||_E\to 0$ as $j\to\infty$ for all $E\subset\subset\O$. Then $u_j\to u$ in $C_n$-capacity on every $E\subset\subset\O$ as $j\to\infty$.}
\medskip
\noindent
{\sl Proof.} Let $\O '\subset\subset\O$ and $\delta >0$. Put
$$A_j=\{z\in\overline{\O '}\ :\ |u_j-u|\geq\delta\}=\{z\in\overline{\O '}\ :\ u-u_j\geq\delta\}.$$
We prove that $C_n(A_j)\to 0$ as $j\to\infty$. Given $\e>0$. By quasicontinuity of $u$ and $u_j$, there is an open set $G$ in $\O$ such that $C_n(G)<\e$, and $u_j|_{\O\backslash G},\ u|_{\O\backslash G}$ are continuous. We have
$$A_j=B_j\cup\{z\in G\ :\ u-u_j\geq\delta\}.$$
where $B_j=\{z\in\overline{\O '}\backslash G\ :\ u-u_j\geq\delta\}$ are compact sets in $\O$ and
$$\varlimsup\limits_{j\to\infty}C_n(A_j)\leq\varlimsup\limits_{j\to\infty}C_n(B_j)+\e$$
We claim that $\lim\limits_{j\to\infty}C_n(B_j)=0$. By Proposition 3.4 we have
$$\aligned C_n(B_j)&=\int\limits_{B_j}(dd^ch_{B_j}^*)^n
\\&\leq \frac 1 {\delta ^n}\int\limits_{B_j}(u-u_j)^n(dd^ch_{B_j}^*)^n
\\&\leq \frac {n!} {\delta ^n} \int\limits_\O(-h_{B_j}^*)[(dd^cu_j)^n-(dd^cu)^n]
\\&\leq \frac {n!} {\delta ^n}\{||(dd^cu_j)^n-(dd^cu)^n||_K+\int\limits_{\O\backslash K}(-h_{\O '})[(dd^cu_j)^n+(dd^cu)^n]\}
\\&\leq \frac {n!} {\delta ^n}\{||(dd^cu_j)^n-(dd^cu)^n||_K+\sup_{\O\backslash K}|h_{\O '}|[\sup\limits_{j\geq 1}\int\limits_\O(dd^cu_j)^n+\int\limits_\O(dd^cu)^n]\}.\endaligned$$
As $\lim\limits_{z\to\partial\O}h_{\O'}(z)=0$ there exists $K\subset\subset\O$ such that
$$\frac {n!} {\delta ^n}\sup_{\O\backslash K}|h_{\O '}|[\sup\limits_{j\geq 1}\int\limits_\O(dd^cu_j)^n+\int\limits_\O(dd^cu)^n]<\e.$$
By the hypothesis
$$\frac {n!} {\delta ^n}||(dd^cu_j)^n-(dd^cu)^n||_K<\e\ \text{for}\ j> j_0.$$
Thus
$$C_n(B_j)<2\e\ \text{for}\ j> j_0.$$
This proves the claim and hence the theorem.
\medskip
\noindent
As an application of Theorem 3.5 we have the following
\medskip
\noindent
{\bf 3.6. Proposition.} {\sl Let $g(A_j)$ be multipolar Green functions on $\O$ such that}
$$\hat A_j =\{w^j_1,...,w^j_{p_j}\}\to\partial\O\ \text{and}\ \sup\limits_{j\geq 1}\nu (A_j)=\sup\limits_{j\geq 1}\sum\limits_{k=1}^{p_j}(\nu^j_k)^n<+\infty$$
{\sl Then $g(A_j)\to 0$ as $j\to\infty$ in $C_n$-capacity.}
\medskip
\noindent
{\sl Proof.} By the hypothesis we have
$$\sup\limits_{j\geq 1}(dd^cg(A_j))^n(\O)=\sup\limits_{j\geq 1}\nu(A_j)<+\infty$$
and
$$||(dd^cg(A_j))^n||_K\to 0\ \text{as}\ j\to\infty\ \text{for all}\ K\subset\subset\O .$$
Theorem 3.5 implies that $g(A_j)\to 0$ as $j\to\infty$ in $C_n$-capacity.
\medskip
\noindent
This section ends up with a criterion for pluripolarity
\medskip
\noindent
{\bf 3.7. Theorem.} {\sl Let $u_j\in\Cal F$ such that $\sup\limits_{j\geq 1}\int\limits_\O (dd^cu_j)^n<+\infty$.}
\medskip
\noindent
{\sl Then there is a constant $A>0$ such that}
\vskip1cm
\medskip
\noindent
i)$(\varlimsup\limits_{j\to\infty} u_j)^*\in\Cal F.$
\medskip
\noindent
ii)$C_n(\{z\in\O :\ (\varlimsup\limits_{j\to\infty} u_j)^*(z)<-t\})\leq \frac {A} {t^n}.$
\medskip
\noindent
iii)$\{z\in\O :\ \lim\limits_{j\to\infty}u_j(z)=-\infty\}$ {\sl is pluripolar.}
\medskip
\noindent
{\sl Proof.} i) For each $j\geq 1$ put $v_j=\sup\{u_j,u_{j+1},...\}.$ By [Ce2] $v_j^*\in\Cal F$ and
$$\sup\limits_{j\geq1 }\int\limits_\O (dd^cv_j^*)^n\leq\sup\limits_{j\geq1 }\int\limits_\O (dd^cu_j)^n<+\infty.$$
By [Ce2] we have $v_j^*\searrow v\in\Cal F$.
\medskip
\noindent
ii) By Proposition 3.1 in [CKZ] we have
$$C_n\{z\in\O :\ (\varlimsup\limits_{j\to\infty} u_j)^*(z)<-t\}=C_n\{z\in\O :\ v(z)<-t\}\leq\frac {2^ne_0(v)} {t^n}=\frac A {t^n},$$
where $A=2^n e_0(v).$
\medskip
\noindent
iii) According to [BT2] we have
$$C_n\{z\in\O :\ \lim\limits_{j\to\infty} u_j(z)=-\infty\} =C_n\{z\in\O :\ v(z)=-\infty\}=0.$$
{\sl Remark.} Theorem 3.7 in the case where $u_j$ are multipole Green functions was proved by D.Coman, N.Levenberg and A.Poletsky in Theorem 4.1 of [CLP].
\vskip1cm
\noindent
{\bf 4. Some properties of the Cegrell's classes and applications}
\medskip
\noindent
In this section, first we prove the following
\medskip
\noindent
{\bf 4.1. Theorem.} {\sl Let} $u,u_1,...,u_{n-1}\in \Cal E$, $v\in\text{PSH}^-(\O)$ {\sl and $T=dd^cu_1\wedge ...\wedge dd^cu_{n-1}$. Then}
$$dd^c\max(u,v)\wedge T|_{\{u>v\}}=dd^cu\wedge T|_{\{u>v\}}.$$
We need the following well-known fact.
\medskip
\noindent
{\bf 4.2. Lemma.} {\sl Let $\mu$ be a measure on $\O$ and $f:\O\to\bold R$ a measurable function on $\O$. The following are equivalent}
\medskip
\noindent
i)$\mu (E)=0$ {\sl for all Borell sets} $E\subset\{f\ne 0\}.$
\medskip
\noindent
ii)$\int\limits_E fd\mu =0$ {\sl for every measurable set $E$ in $\O$.}
\medskip
\noindent
{\sl Proof}. i)$\Rightarrow$ii) follows from:
$$\int\limits_E fd\mu =\int\limits_{E\backslash\{f=0\}} fd\mu +\int\limits_{E\cap\{f=0\}}fd\mu =0$$
ii)$\Rightarrow$i). It suffices to show that $\mu=0$ on every $X_\delta =\{f>\delta>0\}$. By the Hahn decomposition theorem, there exist measurable subsets $X_\delta^+$ and $X_\delta^-$ of $X_\delta$ such that $X_\delta = X_\delta^+\cup X_\delta^-$, $X_\delta^+\cap X_\delta^- =\emptyset$ and $\mu\geq 0$ on $X_\delta^+$, $\mu\leq 0$ on $X_\delta^-$. We have
$$\cases \delta\mu(X_\delta^+)\leq\int\limits_{X_\delta^+}fd\mu =0\\ \delta\mu(X_\delta^-)\geq\int\limits_{X_\delta^-}fd\mu =0\endcases$$
Hence, $\mu(X_\delta^+)=\mu(X_\delta^-)=0$. Therefore, we have $\mu =0$ on $X_\delta$.
\medskip
\noindent
{\sl Proof of Theorem 4.1.}
\medskip
\noindent
a) First we prove the proposition for $v\equiv a<0$. According to the remark following Definition 4.6 in [Ce2], without loss of generality we may assume that $u,u_1,...,u_{n-1}\in\Cal F$. Using Theorem 2.1 in [Ce2] we can find
$$\Cal E_0\cap C(\bar\O)\ni u^j\searrow u,\ \Cal E_0\cap C(\bar\O)\ni u_k^j\searrow u_k,\ k=1,...,n-1.$$
\medskip
\noindent
Since $\{u^j>a\}$ is open we have
$$dd^c\max(u^j,a)\wedge T_j|_{\{u^j>a\}}=dd^cu^j\wedge T_j|_{\{u^j>a\}}.$$
Thus from the inclusion $\{u>a\}\subset\{u^j>a\}$ we obtain
$$dd^c\max(u^j,a)\wedge T_j|_{\{u>a\}}=dd^cu^j\wedge T_j|_{\{u>a\}}.$$
where $T_j=dd^cu_1^j\wedge ...\wedge dd^cu_{n-1}^j.$ By Corollary 5.2 in [Ce2], it follows that
$$\max(u-a,0)dd^c\max(u^j,a)\wedge T_j\to \max(u-a,0)dd^c\max(u,a)\wedge T.$$
$$\max(u-a,0)dd^cu^j\wedge T_j\to \max(u-a,0)dd^cu\wedge T.$$
Hence
$$\max(u-a,0)[dd^c\max(u,a)\wedge T-dd^cu\wedge T]=0.$$
Using Lemma 4.2 we have
$$dd^c\max(u,a)\wedge T=dd^cu\wedge T\ \text{on}\ \{u>a\}.$$
b) Assume that $v\in PSH^-(\O)$. Since $\{u>v\}=\bigcup\limits_{a\in\bold Q^-}\{u>a>v\}$, it suffices to show that
$$dd^c\max(u,v)\wedge T=dd^cu\wedge T\ \text{on}\ \{u>a>v\}$$
for all $a\in\bold Q^-$.
Since $\max(u,v)\in\Cal E$, by a) we have
\medskip
\noindent
$$\aligned dd^c\max(u,v)\wedge T|_{\{\max(u,v)>a\}}&=dd^c\max(\max(u,v),a)\wedge T|_{\{\max(u,v)>a\}}\\&=dd^c\max(u,v,a)\wedge T|_{\{\max(u,v)>a\}}.\endaligned\tag 1$$
$$dd^cu\wedge T|_{\{u>a\}}=dd^c\max(u,a)\wedge T|_{\{u>a\}}.\tag 2$$
Since $\max (u,v,a)=\max (u,a)$ on set open $\{a>v\}$ , we have
$$dd^c\max(u,v,a)\wedge T|_{\{a>v\}}=dd^c\max(u,a)\wedge T|_{\{a>v\}}.\tag 3$$
Since $\{u>a>v\}\subset\{u>a\},\ \{a>v\},\ \{\max(u,v)>a\}$ and (1), (2), (3) we have
$$dd^c\max(u,v)\wedge T|_{\{u>a>v\}}=dd^cu\wedge T|_{\{u>a>v\}}.$$
The next result is an analogue of an inequality due to Demaily in [De2]
\medskip
\noindent
{\bf 4.3. Proposition.} a) {\sl $u,v\in\Cal E$ such that $(dd^cu)^n(\{u=v=-\infty\})=0$. Then}
$$(dd^c\max(u,v))^n\geq 1_{\{u\geq v\}}(dd^cu)^n+1_{\{u<v\}}(dd^cv)^n$$
{\sl where $1_E$ denotes the characteristic function of $E$.}
\medskip
\noindent
b) {\sl Let $\mu$ be a positive measure which vanishes on all pluripolar subsets of $\O$. Suppose $u,v\in\Cal E$ such that $(dd^cu)^n\geq\mu,(dd^cv)^n\geq\mu$. Then $(dd^c\max (u,v))^n\geq\mu$.}
\medskip
\noindent
{\sl Proof.} a) For each $\epsilon >0$ put $A_\epsilon=\{u=v-\epsilon\}\backslash\{u=v=-\infty\}$. Since $A_\epsilon\cap A_\delta =\emptyset$ for $\epsilon\ne\delta$ there exists $\epsilon_j \searrow 0$ such that $(dd^cu)^n(A_{\e_j})=0$ for $j\geq 1$. On the other hand, since $(dd^cu)^n(\{u=v=-\infty\})=0$ we have $(dd^cu)^n(\{u=v-\e_j\})=0$ for $j\geq 1$. Since Theorem 4.1 it follows that
$$\aligned (dd^c\max(u,v-\e_j))^n&\geq (dd^c\max(u,v-\e_j))^n|_{\{u>v-\e_j\}} +(dd^c\max(u,v-\e_j))^n|_{\{u<v-\e_j\}}
\\&=(dd^cu)^n|_{\{u\geq v-\e_j\}}+(dd^cv)^n|_{\{u<v-\e_j\}}
\\&=1_{\{u\geq v-\e_j\}}(dd^cu)^n+1_{\{u< v-\e_j\}}(dd^cv)^n
\\&\geq 1_{\{u\geq v\}}(dd^cu)^n+1_{\{u< v-\e_j\}}(dd^cv)^n.\endaligned$$
Letting $j\to\infty$ and by Remark under Theorem 5.15 in [Ce2] we get
$$(dd^c\max(u,v))^n\geq 1_{\{u\geq v\}}(dd^cu)^n+1_{\{u<v\}}(dd^cv)^n$$
because $\max(u,v-\e_j)\nearrow\max (u,v)$ and $1_{\{u<v-\e_j\}}\nearrow 1_{\{u<v\}}$ as $j\to\infty$.
\medskip
\noindent
b) Argument as a)
\medskip
\noindent
{\bf 4.4. Proposition.} {\sl Let} $u_1,...,u_k\in\text{PSH}(\O)\cap L^\infty(\O)$ {\sl and $u_{k+1},...,u_n\in\Cal E$. Then}
\medskip
\noindent
i) $\int\limits_Bdd^cu_1\w ...\w dd^cu_n=O((C_n(B))^{\frac k n})$ {\sl for all Borel sets } $B\subset\O'\subset\subset\O.$
\medskip
\noindent
ii) $\int\limits_{B(a,r)}dd^cu_1\w ...\w dd^cu_n=o((C_n(B(a,r)))^{\frac k n})$ {\sl as} $r\to 0$ {\sl for all} $a\in\O.$
\medskip
\noindent
where $B(a,r)=\{z\i\bold C^n\ :\ |z-a|<r\}$
\medskip
\noindent
{\sl Proof.} We may assume that $0\leq u_j\leq 1$ for $j=1,...,k$. On the other hand, by the remark following Defintion 4.6 in [Ce2] we again may assume that $u_{k+1},...,u_n\in\Cal F.$
\medskip
\noindent
i) For each open set $B\subset\subset\O$, applying Proposition 3.1 we get
$$\aligned\int\limits_Bdd^cu_1\w ...\w dd^cu_n&=\int\limits_B(-h_B^*)^kdd^cu_1\w ...\w dd^cu_n
\\&\leq\int\limits_\O(-h_B^*)^kdd^cu_1\w ...\w dd^cu_n
\\&\leq k!\int\limits_\O(1-u_1)(dd^ch_B^*)^k\w dd^cu_{k+1}\w ...\w dd^cu_n
\\&\leq k!\int\limits_\O(dd^ch_B^*)^k\w dd^cu_{k+1}\w ...\w dd^cu_n
 \\&\leq k![\int\limits_\O(dd^ch_B^*)^n]^{\frac k n}\w [\int\limits_\O(dd^cu_{k+1})^n]^{\frac 1 n}\w ...\w [\int\limits_\O(dd^cu_n)^n]^{\frac 1 n}
\\&\text{(by Corollary 5.6 in [Ce2])}
\\&\leq k!(e_0(u_{k+1}))^{\frac 1 n}...(e_0(u_n))^{\frac 1 n}.[C_n(B)]^{\frac k n}
\\&\leq\text{constants}.[C_n(B)]^{\frac k n}.\endaligned$$
Hence
$$\int\limits_Bdd^cu_1\w ...\w dd^cu_n\leq\text{constants}.[C_n(B)]^{\frac k n}.$$
for all Borel set $B\subset\O$.
\medskip
\noindent
ii) By Proposition 3.1 we have
$$\aligned\int\limits_\O(-\varphi)^kdd^cu_1\w ...\w dd^cu_n&\leq k!\int\limits_\O(1-u_1)(dd^c\varphi)^k\w dd^cu_{k+1}\w ...\w dd^cu_n
\\&\leq k!\int\limits_\O(dd^c\varphi)^k\w dd^cu_{k+1}\w ...\w dd^cu_n<+\infty.\endaligned$$
\medskip
\noindent
Hence $(-\varphi)^k\in L_1(dd^cu_1\w ...\w dd^cu_n)$ for all $\varphi\in\Cal F(\O)$. Given $a\in\O$ let $r_0,R_0$ such that $B(a,r_0)\subset\subset\O\subset\subset B(a,R_0)$. Then
$$\log\frac {|z-a|} {R_0}\leq g_a(z)\leq\log\frac {|z-a|} {r_0}$$
for all $z\in\O$, where $g_a$ denotes the Green function of $\O$ with pole at $a$. Since $(-g_a)^k\in L_1(dd^cu_1\w ...\w dd^cu_n)$, it follows that
$$\int\limits_{B(a,r)}(-g_a)^kdd^cu_1\w ...\w dd^cu_n\to 0\ \text{as}\ r\to 0$$
Hence
$$(\log r_0-\log r)^k\int\limits_{B(a,r)}dd^cu_1\w ...\w dd^cu_n\leq\int\limits_{B(a,r)}(-g_a)^kdd^cu_1\w ...\w dd^cu_n\to 0$$
as $r\to 0$. This means that
$$\int\limits_{B(a,r)}dd^cu_1\w ...\w dd^cu_n=o((\frac 1 {\log r_0-\log r})^k)\ \text{as}\ r\to 0$$
Combining this with the inequality
$$C_n(B(a,r),\O)\geq C_n(B(a,r),B(a,R_0))=(\frac {1} {\log R_0-\log r})^n=O((\frac {1} {\log r_0-\log r)^n})$$
we get
$$\int\limits_{B(a,r)}dd^cu_1\w ...\w dd^cu_n=o((C_n(B(a,r)))^{\frac k n}).$$
The next result should be compared with Theorem 6.3 in [Ce1]
\medskip
\noindent
{\bf 4.5. Theorem.} {\sl Let $u_1,...,u_n\in\Cal E$. Then there exists $\tilde u\in\Cal E^a$ such that}
$$dd^cu_1\w ...\w dd^cu_n=(dd^c\tilde u)^n+dd^cu_1\w ...\w dd^cu_n|_{\{u_1=...=u_n=-\infty\}}.$$
{\sl Proof.} First, we write
$$dd^cu_1\w ...\w dd^cu_n=\mu +dd^cu_1\w ...\w dd^cu_n|_{\{u_1=...=u_n=-\infty\}}.$$
where
$$\mu=dd^cu_1\w ...\w dd^cu_n|_{\{u_1>-\infty\}\cup ...\cup\{u_n>-\infty\}}.$$
It is easy to see that $\mu\ll C_n$ in every $E\subset\subset\O$. Indeed, by Theorem 4.1 we have
$$dd^cu_1\w ...\w dd^cu_n|_{\{u_1>-j\}}=dd^c\max(u_1,-j)\w ...\w dd^cu_n|_{\{u_1>-j\}}.$$
Hence, by Proposition 4.4 (i) it follows that $dd^cu_1\w ...\w dd^cu_n|_{\{u_1>-j\}}\ll C_n$ in every $E\subset\subset\O$. Next, it remains to show that there exists $\tilde u\in\Cal E^a$ such that $\mu=(dd^c\tilde u)^n$. Let $\{\O_j\}$ be an increasing exhaustion sequence of $\O$. For each $j\geq 1$ put $\mu_j=\mu|_{\O_j}$. By [\AA h] there exists $\tilde u_j\in\Cal F$ such that $(dd^c\tilde u_j)^n=\mu_j$. Notice that $\mu_j\nearrow\mu$ and
$$(dd^c\tilde u_j)^n\leq \mu\leq (dd^c(u_1+...+u_n))^n.$$
Applying the comparison principle we obtain
$$\tilde u_j\searrow\tilde u\geq u_1+...+u_n\in\Cal E.$$
Hence, $\tilde u\in\Cal E^a$ and $(dd^c\tilde u)^n=\lim\limits_{j\to\infty}(dd^c\tilde u_j)^n=\mu$. The proof is thereby completed.
\medskip
\noindent
{\bf 4.6. Corollary.} {\sl $u_1,...,u_n\in\Cal E$. Then the following are equivalent}
\medskip
\noindent
i) $dd^cu_1\w ...\w dd^cu_n\ll C_n$ {\sl in every} $E\subset\subset\O.$
\medskip
\noindent
ii) $\int\limits_{\{u_1=...=u_n=-\infty\}}dd^cu_1\w...\w dd^cu_n=0.$
\medskip
\noindent
iii) $\int\limits_{\{u_1<-s,...,u_n<-s\}\cap E}dd^cu_1\w...\w dd^cu_n\to 0$ {\sl as} $s\to+\infty$ {\sl for all} $E\subset\subset\O$.
\medskip
\noindent
{\sl Proof.} Direct application of Theorem 4.5.
\medskip
\noindent
The comparison principle for class $\Cal F$ was studied in [Ce3] and [H1]. By using Proposition 3.1 and Theorem 4.1 we prove a Xing type comparison principle for $\Cal F$
\medskip
\noindent
{\bf 4.7. Theorem.} {\sl Let $u\in\Cal F$, $v\in\Cal E$ and $1\leq k\leq n$. Then}
$$\frac 1 {k!} \int\limits_{\{u<v\}}(v-u)^kdd^cw_1\w ...\w dd^cw_n+\int\limits_{\{u<v\}}(r-w_1)(dd^cv)^k\w dd^cw_{k+1}\w ...\w dd^cw_n$$
$$\leq \int\limits_{\{u<v\}\cup\{u=v=-\infty\}}(r-w_1)(dd^cu)^k\w dd^cw_{k+1}\w ...\w dd^cw_n$$
{\sl for all} $w_j\in\text{PSH}(\O),\ 0\leq w_j\leq 1,\ j=1,...,k,\ w_{k+1},...,w_n\in\Cal F$ {\sl and all $r\geq 1$.}
\medskip
\noindent
{\sl Proof.} Let $\e>0$. We set $\tilde v=\max (u,v-\e)$. By a) in Proposition 3.1 we have
$$\frac 1 {k!} \int\limits_\O(\tilde v-u)^kdd^cw_1\w ...\w dd^cw_n+\int\limits_\O(r-w_1)(dd^c\tilde v)^k\w dd^cw_{k+1}\w ...\w dd^cw_n$$
$$\leq \int\limits_\O(r-w_1)(dd^cu)^k\w dd^cw_{k+1}\w ...\w dd^cw_n.$$
Since $\{u< \tilde v\}=\{u<v-\e\}$ and Theorem 4.1  we have
$$\frac 1 {k!} \int\limits_{\{u<v-\e\}}(v-\e-u)^kdd^cw_1\w ...\w dd^cw_n+\int\limits_{\{u\leq v-\e\}}(r-w_1)(dd^cv)^k\w dd^cw_{k+1}\w ...\w dd^cw_n$$
$$\leq \int\limits_{\{u\leq v-\e\}}(r-w_1)(dd^cu)^k\w dd^cw_{k+1}\w ...\w dd^cw_n$$
$$\leq \int\limits_{\{u<v\}\cup\{u=v=-\infty\}}(r-w_1)(dd^cu)^k\w dd^cw_{k+1}\w ...\w dd^cw_n.$$
Letting $\e\searrow 0$ we obtain
$$\frac 1 {k!} \int\limits_{\{u<v\}}(v-u)^kdd^cw_1\w ...\w dd^cw_n+\int\limits_{\{u<v\}}(r-w_1)(dd^cv)^k\w dd^cw_{k+1}\w ...\w dd^cw_n$$
$$\leq \int\limits_{\{u<v\}\cup\{u=v=-\infty\}}(r-w_1)(dd^cu)^k\w dd^cw_{k+1}\w ...\w dd^cw_n$$
{\bf 4.8. Corollary.} {\sl Let $u\in\Cal E^a$ such that $u\geq v$ for all functions $v\in\Cal E$ satisfying $(dd^cu)^n\leq (dd^cv)^n$. Then}
\medskip
\noindent
 $$\frac 1 {n!} \int\limits_{\{u<v\}}(v-u)^ndd^cw_1\w ...\w dd^cw_n+\int\limits_{\{u<v\}}(r-w_1)(dd^cv)^n$$
\medskip
\noindent
$$\leq\int\limits_{\{u<v\}}(r-w_1)(dd^cu)^n$$
\medskip
\noindent
{\sl for all $v\in\Cal E$, $r\geq 1$ and all} $w_1,...,w_n\in\text{PSH}(\O)$, $0\leq w_1,...,w_n\leq 1$.
\medskip
\noindent
{\sl Proof.} Let $\{\O_j\}$ be an increasing exhaustion sequence of relatively compact subdomains of $\O$. Set $\mu_j=1_{\O_j}1_{\{u>-j\}}(dd^cu)^n$, where $1_E$ denotes the characteristic function of $E\subset\O$. Applying Theorem 4.1 we have
$$\mu_j=1_{\O_j}1_{\{u>-j\}}(dd^c\max(u,-j))^n\leq 1_{\O_j}(dd^c\max(u,-j))^n.$$
Take $\phi\in\Cal E_0(\O)\cap C(\bar\O)$. Put
$$\phi_j=\max(u,-j,a_j\phi)$$
where $a_j=\frac {-j} {\sup\limits_{\O_{j+1}}\phi}.$ Then $\phi_j=\max(u,-j)$ on $\O_{j+1}$, $\phi_j\in\Cal E_0$ and
$$\mu_j\leq 1_{\O_j}(dd^c\max(u,-j))^n=1_{\O_j}(dd^c\phi_j)^n\leq (dd^c\phi_j)^n.$$
By Ko{\l}odziej's theorem (see [Ko]) there exists $u_j\in\Cal E_0$ such that
$$(dd^cu_j)^n=\mu_j=1_{\O_j}1_{\{u>-j\}}(dd^cu)^n, \forall\ j\geq 1.$$
for all $j\geq 1$. By the comparison principle we have $u_j\searrow\tilde u\geq u$. On the other hand, since $(dd^cu)^n(\{u=-\infty\})=0$, it follows that
$$(dd^cu_j)^n=1_{\O_j}1_{\{u>-j\}}(dd^cu)^n\to (dd^cu)^n$$
weakly as $j\to\infty.$ Thus $(dd^c\tilde u)^n=\lim\limits_{j\to\infty}(dd^cu_j)^n=(dd^cu)^n$. By the hypothesis we have $\tilde u=u$. Applying Theorem 4.7 we get
$$\frac 1 {n!} \int\limits_{\{u_j<v\}}(v-u_j)^ndd^cw_1\w ...\w dd^cw_n+\int\limits_{\{u_j<v\}}(r-w_1)(dd^cv)^n$$
$$\leq\int\limits_{\{u_j<v\}}(r-w_1)(dd^cu_j)^n$$
$$\leq\int\limits_{\{u_j<v\}}(r-w_1)(dd^cu)^n.$$
Letting $j\to\infty$ we obtain
$$\frac 1 {n!} \int\limits_{\{u<v\}}(v-u)^ndd^cw_1\w ...\w dd^cw_n+\int\limits_{\{u<v\}}(r-w_1)(dd^cv)^n$$
Arguing as in Theorem  4.7 we prove a Xing type comparison principle for $\Cal E$.
\medskip
\noindent
{\bf 4.9. Theorem.} {\sl Let $u,v\in\Cal E$ and $1\leq k\leq n$ such that $\varliminf\limits_{z\to\partial\O}[u(z)-v(z)]\geq 0$. Then}
$$\frac 1 {k!} \int\limits_{\{u<v\}}(v-u)^kdd^cw_1\w ...\w dd^cw_n+\int\limits_{\{u<v\}}(r-w_1)(dd^cv)^k\w dd^cw_{k+1}\w ...\w dd^cw_n$$
$$\leq \int\limits_{\{u<v\}\cup\{u=v=-\infty\}}(r-w_1)(dd^cu)^k\w dd^cw_{k+1}\w ...\w dd^cw_n$$
{\sl for all} $w_j\in\text{PSH}(\O),\ 0\leq w_j\leq 1,\ j=1,...,k,\ w_{k+1},...,w_n\in\Cal E$ {\sl and all $r\geq 1$.}
\medskip
\noindent
{\sl Proof.} Let $\e>0$. We set $\tilde v=\max (u,v-\e)$. By b) in Proposition 3.1 we have
$$\frac 1 {k!} \int\limits_\O(\tilde v-u)^kdd^cw_1\w ...\w dd^cw_n+\int\limits_\O(r-w_1)(dd^c\tilde v)^k\w dd^cw_{k+1}\w ...\w dd^cw_n$$
$$\leq \int\limits_\O(r-w_1)(dd^cu)^k\w dd^cw_{k+1}\w ...\w dd^cw_n.$$
Since $\{u< \tilde v\}=\{u<v-\e\}$ and Theorem 4.1  we have
$$\frac 1 {k!} \int\limits_{\{u<v-\e\}}(v-\e-u)^kdd^cw_1\w ...\w dd^cw_n+\int\limits_{\{u\leq v-\e\}}(r-w_1)(dd^cv)^k\w dd^cw_{k+1}\w ...\w dd^cw_n$$
$$\leq \int\limits_{\{u\leq v-\e\}}(r-w_1)(dd^cu)^k\w dd^cw_{k+1}\w ...\w dd^cw_n$$
$$\leq \int\limits_{\{u<v\}\cup\{u=v=-\infty\}}(r-w_1)(dd^cu)^k\w dd^cw_{k+1}\w ...\w dd^cw_n.$$
Letting $\e\searrow 0$ we obtain
$$\frac 1 {k!} \int\limits_{\{u<v\}}(v-u)^kdd^cw_1\w ...\w dd^cw_n+\int\limits_{\{u<v\}}(r-w_1)(dd^cv)^k\w dd^cw_{k+1}\w ...\w dd^cw_n$$
$$\leq \int\limits_{\{u<v\}\cup\{u=v=-\infty\}}(r-w_1)(dd^cu)^k\w dd^cw_{k+1}\w ...\w dd^cw_n.$$
\vskip1cm
\noindent
\centerline{\bf References}
\vskip1cm
\noindent
[\AA h] P. \AA hag, {\sl The complex Monge-Amp\`{e}re operator on bounded hyperconvex domains}, Ph. D. Thesis, Ume\aa\ University, (2002).
\medskip
\noindent
[Bl1] Z. Blocki, {\sl On the definition of the Monge-Amp\`{e}re operator in $\bold C^2$}, Math. Ann., {\bf 328} (2004), 415-423.
\medskip
\noindent
[Bl2] Z. Blocki, Weak solutions to the complex Hessian equation, Ann. Inst. Fourier {\bf 55} (2005), 1735-1756.
\medskip
\noindent
[BT1] E. Bedford and B.A.Taylor, {\sl The Dirichlet problem for the complex Monge-Amp\`{e}re operator}. Invent. Math.{\bf 37} (1976), 1-44.
\medskip
\noindent
[BT2] E. Bedford and B.A.Taylor, {\sl A new capacity for plurisubharmonic functions}. Acta Math., {\bf 149} (1982), 1-40.
\medskip
\noindent
[BT3] E. Bedford and B.A.Taylor, {\sl Fine topology, Silov boundary, and $(dd^c)^n$}. J. Funct. Anal. {\bf 72} (1987), 225-251.
\medskip
\noindent
[Ce1] U. Cegrell, {\sl Pluricomplex energy}. Acta Math., {\bf 180} (1998), 187-217.
\medskip
\noindent
[Ce2] U. Cegrell, {\sl The general definition of the complex Monge-Amp\`{e}re operator}. Ann. Inst. Fourier (Grenoble) {\bf 54} (2004), 159-179.
\medskip
\noindent
[Ce3] U. Cegrell, {\sl A general Dirichlet problem for the complex Monge-Amp\`{e}re operator}, preprint (2006).
\medskip
\noindent
[CKZ] U. Cegrell, S. Ko{\l}odziej and A. Zeriahi, {\sl Subextention of plurisubharmonic functions with weak singularities}. Math. Zeit., {\bf 250} (2005), 7-22.
\medskip
\noindent
[Cz] R. Czyz, {\sl Convergence in capacity of the Perron-Bremermann envelope}, Michigan Math. J., {\bf 53} (2005), 497-509.
\medskip
\noindent
[CLP] D. Coman, N. Levenberg and E.A. Poletsky,{\sl Quasianalyticity and pluripolarity}, J. Amer. Math. Soc., {\bf 18} (2005), 239-252.
\medskip
\noindent
[De1] J-P. Demailly, {\sl Monge-Amp\`{e}re operators, Lelong Numbers and Intersection theory}, Complex Analysis and Geometry, Univ. Ser. Math., Plenum, New York, 1993, 115-193.
\medskip
\noindent
[De2] J-P. Demailly, {\sl Potential theory in several variables}, preprint (1989).
\medskip
\noindent
[Ko] S. Ko{\l}odziej, {\sl The range of the complex Monge-Amp\`{e}re operator}, II, Indiana Univ. Math. J., {\bf 44} (1995), 765-782.
\medskip
\noindent
[H1] P. Hiep, {\sl  A characterization of bounded plurisubharmonic functions}, Ann. Polon. Math., {\bf 85} (2004), 233-238.
\medskip
\noindent
[H2] P. Hiep, {\sl The comparison principle and Dirichlet problem in the class $\Cal E_p(f)$, $p>0$}, Ann. Polon. Math., {\bf 88} (2006), 247-261.
\medskip
\noindent
[Le] P .Lelong, {\sl Notions capacitaires et fonctions de Green pluricomplexes dans les espaces de Banach}. C.R. Acad. Sci. Paris Ser. Imath., 305:71-76, 1987.
\medskip
\noindent
[Xi1] Y. Xing, {\sl Continuity of the complex Monge-Amp\`{e}re operator}. Proc. of Amer. Math. Soc., {\bf 124} (1996), 457-467.
\medskip
\noindent
[Xi2] Y. Xing, {\sl Complex Monge-Amp\`{e}re measures of pluriharmonic functions with bounded values near the boundary}. Cand. J. Math., {\bf 52}, (2000),1085-1100.
\vskip1cm
\noindent
Department of Mathematics
\medskip
\noindent
Hanoi University of Education (Dai hoc Su Pham Hanoi).
\medskip
\noindent
Cau giay, Ha Noi, VietNam
\medskip
\noindent
E-mail: phhiep$_-$vn$\@$yahoo.com
\end